\documentclass[11pt,twoside,a4paper]{article}
 \usepackage{euscript,a4,times}
  \usepackage[english]{babel}
   \usepackage{amsfonts,amssymb,amsmath,amsthm}
    \usepackage{latexsym,mathrsfs,color}
\catcode`@=11
\newbox\tr@tto
\setbox\tr@tto=\hbox{{\count0=0\dimen0=-,9pt\dimen1=1,1pt\loop\ifnum\count0<11
\advance \count0 by1 \vrule width.51pt height\dimen1
     depth\dimen0\kern-0.17pt\advance\dimen0 by-0.05pt\advance\dimen1
     by0.1pt\repeat \loop\ifnum\count0<21\advance \count0 by1 \vrule
     width.6pt height\dimen1 depth\dimen0\kern-0.2pt \advance\dimen0
     by-0.1pt\advance\dimen1 by 0.05pt\repeat}}
\def\medint{\displaystyle\copy\tr@tto\kern-10.4pt\int}
\catcode`@=12
\input amssym.def
\newcommand{\R}{\mathbb{R}}
 \newcommand{\Rn}{\mathbb{R}^n}
  \newcommand{\N}{\mathbb{N}}
   \newcommand{\Om}{\Omega}
    \newcommand{\test}{C_c^1}
     \newcommand{\prt}{\partial}
      \newcommand{\omn}{\omega_{n}}
       \newcommand{\dd}{\mathrm{d}}
        \newcommand{\Ln}{\mathcal{L}^n}
         
          \newcommand{\e}{\varepsilon}

             \newcommand{\Div}{\operatorname{div}}

                    \newcommand{\ub}{{\mathbb{B}^n}}

                          \newcommand{\dS}{\dd\sigma}
                           \newcommand{\us}{{\mathbb{S}^{n-1}}}
                             
\def\aver{\mathop{\int\mkern-18.5mu{-}}\nolimits}
\newtheorem{theorem}{Theorem}[section]
 \newtheorem{proposition}[theorem]{Proposition}
  \newtheorem{lemma}[theorem]{Lemma}
   
    \theoremstyle{definition}

       \newtheorem{remark}[theorem]{Remark}
        
         \newtheorem{acknowledgements}{Acknowledgements}
\begin{document}
\title{\bf{\Large{The optimal Leray-Trudinger inequality}}}
\author{\small{GIUSEPPINA DI BLASIO}$^\ast$\\ \small{giuseppina.diblasio@unicampania.it}\\\\ \small{GIOVANNI PISANTE}$^\ast$\\ \small{giovanni.pisante@unicampania.it}\\\\ \small{GEORGIOS PSARADAKIS}$^\ast$\\ \small{georgios.psaradakis@unicampania.it}
\\\\ $^\ast$Dipartimento di Matematica e Fisica\\
	Universit\`a degli Studi della Campania ``L. Vanvitelli"\\
	Viale Lincoln 5, 81100 Caserta, Italy}
\maketitle

\tableofcontents
\begin{abstract}
\noindent We fill the gap left open in \cite{MT}, regarding the minimum exponent on the logarithmic correction weight so that the Leray-Trudinger inequality (see \cite{PsSp}) holds. Instead of the representation formula used in \cite{PsSp} and \cite{MT}, our proof uses expansion in spherical harmonics as in \cite{VzZ}.
\end{abstract}

\section{Introduction}

Let $\Om$ be an open subset of $\Rn$, $n\in\N\setminus\{1,2\}$. The classical Sobolev inequality asserts that
\begin{align}\nonumber
 \sup_{u\in\mathcal{D}_n(\Om)}\int_{\Om}|u|^{2n/(n-2)}~\dd x
  <\infty,
\end{align}
where we have set
\begin{align}\nonumber
 \mathcal{D}_n(\Om):=\bigg\{u\in\test(\Om)~\Big|~\int_\Om|\nabla u|^2~\dd x\leq1\bigg\}.
\end{align}
It is well known that the exponent $2n/(n-2)$ cannot be increased. At least when $\Om$ has finite Lebesgue measure, this may suggest that functions in $\mathcal{D}_2(\Om)$ could be bounded. Standard examples show that this is not the case and Neil Trudinger in \cite{Tr} has established the optimal embedding in this case. More precisely, \textit{Trudinger's inequality} (see \cite{Pe}, \cite{Po} and \cite{Y} for prior results and \cite{M} for the best constant) says that there exists a positive constant $c$ such that
\begin{align}\label{Trudinger inequality 2d}
 \sup_{u\in\mathcal{D}_2(\Om)}
  \aver_{\hspace{-.3em}\Om}e^{\alpha u^2}~\dd x
  <\infty\qquad\forall~\alpha<c,
\end{align}
and the exponent in the power $u^2$ cannot be increased\footnote{In \eqref{Trudinger inequality 2d} and throughout this paper $\aver_\Om f~\dd x$ stands for $\big(\Ln(\Om)\big)^{-1}\int_\Om f~\dd x$.}.

\smallskip

Consider now the following higher dimensional Hardy-type inequality
\begin{align}\label{Hardy inequality}
 \int_\Om|\nabla u|^2~\dd x
  -\Big(\frac{n-2}{n}\Big)^2\int_\Om\frac{|u|^2}{|x|^2}~\dd x\geq0
   \qquad\forall~u\in\test(\Om).
\end{align}
Leray seems to be the first to have provided a proof of \eqref{Hardy inequality} for $n=3$ (see \cite[pg 47]{Le}). Note that in case $0\in\Om$, the constant $\big((n-2)/n\big)^2$ turns out to be the best possible. Nevertheless, already in his paper, Leray proved a substitute inequality for the case $n=2$ (see \cite[pg 49]{Le}). Of interest in this paper is the following version of \textit{Leray's inequality}: if $\Om\subset\R^2$ is a bounded domain that contains the origin, then we have
\begin{align}\label{Leray inequality}
 I_2[u]:=\int_\Om|\nabla u|^2~\dd x
  -\frac14\int_\Om\frac{|u|^2}{|x|^2}X^2_1\Big(\frac{|x|}{R_\Om}\Big)~\dd x
  \geq0
   \qquad\forall~u\in\test(\Om),
\end{align}
where $R_\Om:=\sup_{x\in\Om}|x|$ and \[X_1(t):=(1-\log t)^{-1},\qquad t\in(0,1], \quad X_1(0):=0.\] Moreover, the constant $1/4$ is the best possible and the power $2$ on $X_1$ cannot be decreased (see for example \cite[Theorems 4.3 \& 5.4]{BFT1} with $k=N$ there).

\smallskip

Now let $u\in\test(\Om)\setminus\{0\}$ where $\Om\subset\R^2$ is as above and suppose further that $I_2[u]\leq1$. Setting
\begin{align}\nonumber
 v
  =
   uX_1^{1/2},
\end{align}
and using integration by parts, one easily obtains (see \cite[Proposition 2.6]{PsSp})
\begin{align}\label{transition twodim}
 \int_\Om|\nabla v|^2X_1^{-1}\Big(\frac{|x|}{R_\Om}\Big)~\dd x
  =
   I_2[u].
\end{align}
Since $X_1(t)\leq1$ for all $t\in[0,1]$, equality \eqref{transition twodim} readily implies $\int_\Om|\nabla v|^2\dd x\leq I_2[u]\leq1$. Hence
\begin{align}\nonumber
 \aver_{\hspace{-.3em}\Om}e^{\alpha u^2X_1(|x|/R_\Om)}\dd x
  \leq
   \sup_{v\in\mathcal{D}_2(\Om)}
    \aver_{\hspace{-.3em}\Om}e^{\alpha v^2}\dd x
    <\infty,
\end{align}
because of Trudinger's inequality \eqref{Trudinger inequality 2d}. Consequently, with $c$ as in \eqref{Trudinger inequality 2d} and \[\mathcal{I}_2(\Om):=\big\{u\in\test(\Om)~\big|~I_2[u]\leq1\big\},\] we have the following combination of \eqref{Leray inequality} and \eqref{Trudinger inequality 2d}
\begin{align}\label{elementary Leray-Moser inequality}
 \sup_{u\in\mathcal{I}_2(\Om)}
  \aver_{\hspace{-.3em}\Om}e^{\alpha u^2X_1(|x|/R_\Om)}~\dd x
  <\infty\qquad\forall~\alpha<c.
\end{align}

In this paper we establish the optimal version of \eqref{elementary Leray-Moser inequality}. More precisely, it has been shown in \cite{PsSp} that estimate \eqref{elementary Leray-Moser inequality} is far from being optimal. In fact, \cite[Theorem 1.1]{PsSp} says that given $\e>0$, there exists a positive constant $c=c(\e)$ such that
\begin{align}\nonumber
 \sup_{u\in\mathcal{I}_2(\Om)}
  \aver_{\hspace{-.3em}\Om}e^{\alpha u^2X^\e_1(|x|/R_\Om)}~\dd x
  <\infty\qquad\forall~\alpha<c,
\end{align}
and that such an estimate fails to hold for all $\alpha>0$ when $\e=0$. Inspired by a result of Calanchi and Ruf \cite{CR}, further understanding on the problem was provided by Mallick and Tintarev in \cite{MT}. They showed that for any $\gamma\geq2$ there exists a positive constant $c$, not depending on $\gamma$, such that
\begin{align}\label{MT}
 \sup_{u\in\mathcal{I}_2(\Om)}
  \aver_{\hspace{-.3em}\Om}e^{\alpha u^2X^\gamma_2(|x|/R_\Om)}~\dd x
  <\infty\qquad\forall~\alpha<c,
\end{align}
where
\[
 X_2(\cdot):=X_1\big(X_1(\cdot)\big).
\]
Moreover, such an estimate fails to hold for all $\alpha>0$ when $\gamma<1$. As observed in \cite{MT}, by \cite[Lemma 5]{CR}, inequality \eqref{MT} is true for $\gamma=1$ when restricted to $\mathcal{I}_2^{\rm rad}$, i.e. radially symmetric functions of $\mathcal{I}_2(\mathbb{B}^2)$ \footnote{$\mathbb{B}^n$ stands for the unit ball of $\Rn$ having center at the origin.}. In \S\ref{sec:radial}, we extend this result to the multi-dimensional case. Furthermore, in \S\ref{sec:non-radial} we take away the radial restriction, proving thus the following optimal \emph{Leray-Trudinger inequality}: 

\begin{theorem}\label{main theorem} Let $\Om\subset\Rn$, $n\in\mathbb{N}\setminus\{1\}$, be a bounded domain that contains the origin and set $R_\Om:=\sup_{x\in\Om}|x|$. There exists a positive constant $c=c(n)$ such that
\begin{align}\label{optimal MT multidim}
\sup_{u\in\mathcal{I}_n(\Om)}
\aver_{\hspace{-.3em}\Om} e^{\alpha\big[|u|X^{1/n}_2(|x|/R_\Om)\big]^{n/(n-1)}}~\dd x
  <\infty\qquad\forall~\alpha<c,
\end{align}
where $\mathcal{I}_n(\Om):=\big\{u\in\test(\Om)~\big|~I_n[u]\leq1\big\}$ with
\begin{align}\label{Leray inequality mutidim}
 I_n[u]:=\int_\Om|\nabla u|^n~\dd x
  -\Big(\frac{n-1}{n}\Big)^n\int_\Om\frac{|u|^n}{|x|^n}X^n_1\Big(\frac{|x|}{R_\Om}\Big)~\dd x.
\end{align}
Moreover, the exponent $1/n$ on $X_2$ cannot be increased.
\end{theorem}

Note that $I_n[u]$ is always nonnegative (an elementary proof based on integration by parts and H\"{o}lder's inequality can be found in \cite[Theorem 4.2]{BFT1} 
or \cite[Theorem 2.1]{PsSp}). For the optimality of the exponent on $X_1$ and of the constant which appear in \eqref{Leray inequality mutidim}, we refer to \cite[Theorem 5-(i)]{BFT1}.

It is well known that the proof of the $n$-dimensional Trudinger's inequality:
\begin{align}\label{optimal MT multidim}
\sup_{u\in\mathcal{D}_n(\Om)}
 \aver_{\hspace{-.3em}\Om} e^{\alpha|u|^{n/(n-1)}}~\dd x
  <\infty\qquad\forall~\alpha<c,
\end{align}
for some positive $c=c(n)$, is based on finding the sharp growth of the $L^q$-norm of $W_0^{1,n}$ functions, as $q\rightarrow\infty$. The key estimate to prove \eqref{optimal MT multidim} is indeed
\begin{align}\label{sharp Trudinger growth}
 \Bigg(\aver_{\hspace{-.3em}\Om}|u|^q~\dd x\Bigg)^{1/q}\leq c(n)~q^{1-1/n}\|\nabla u\|_{L^n(\Om)}\qquad\forall~q>n,
\end{align}
whenever $u\in\test(\Om)$. Following a similar path, to establish Theorem \ref{main theorem}, we analogously prove
\begin{align}\label{sharp Leray-Trudinger growth}
 \Bigg(\aver_{\hspace{-.3em}\Om}\Big(|u|X_2^{1/n}(|x|/R_\Om)\Big)^q~\dd x\Bigg)^{1/q}\leq c(n)~q^{1-1/n}\big(I_n[u]\big)^{1/n}\qquad\forall~q>n,
\end{align}
whenever $u\in\test(\Om)$.

For partial results on the corresponding problem dealing with the Hardy inequality that involves the distance to the boundary of convex or mean convex domains, we refer to \cite{WY}, \cite{FP} and \cite{dBPP}.

\section{Preliminary estimates}
\subsection{Lower estimates on $I_n[u;\Om]$}
\textbf{Notation.} From now on we write $X_1$, $X_2$ instead of $X_1(|x|/R_\Om)$, $X_2(|x|/R_\Om)$.

\bigskip

We recall a known lower estimate for the Hardy-Leray difference

\begin{proposition}[\cite{PsSp}-Proposition 2.6]\label{proposition:link:PsSp} Set $\lambda_n:=2^{n-1}-1$. Whenever $u\in\test(\Om)$ we have
\begin{align}\label{link}
 \int_\Om|\nabla v|^nX_1^{-n+1}~\dd x
  \leq
   \lambda_nI_n[u],
\end{align}
where $v:=X_1^{1-1/n}u$. For $n=2$ we have equality in \eqref{link}.
\end{proposition}
In order to prove theorem \ref{main theorem}, we are now going to establish one more lower estimate on $I_n[u;\Om]$. This estimate (see \eqref{link2} below) will be enough to prove our main theorem for radial functions. However, in the next section, it will be combined with Proposition \ref{proposition:link:PsSp} to remove the radiality assumption; see \cite{GkPs}. Observe that for $n=2$, estimate \eqref{link2} agrees with \eqref{link}.

\begin{proposition}\label{proposition:link:Sobolev} Set $\kappa_n:=\lambda_n \big(\frac{2n}{n-1}\big)^{n-2}$. Whenever $u\in\test(\Om)$ we have
\begin{align}\label{link2}
 \int_\Om|x|^{2-n}|v|^{n-2}|\nabla v|^2X_1^{-1}~\dd x
  \leq
   \kappa_nI_n[u],
\end{align}
where $v:=X_1^{1-1/n}u$. For $n=2$ we have equality in \eqref{link2}.
\end{proposition}

\noindent\textbf{Proof.} It suffices to consider $u\in\test\big(\Om\setminus\{0\}\big)$. Setting $u = X_1^{-1+1/n}v$ we compute
\begin{align}\nonumber
 |\nabla u|^n
  =
   \bigg|X_1^{-1+1/n}\nabla v
    -\frac{n-1}{n}X_1^{1/n}\frac{v}{|x|}\frac{x}{|x|}\bigg|^n.
\end{align}
Applying the vectorial inequality (see also \cite{Li})
\begin{align}\label{vec}
 |b-a|^n-|a|^n
  \geq
   \frac{1}{\lambda_n2^{n-2}}|a|^{n-2}|b|^2-n|a|^{n-2}a\cdot b,
\end{align}
we get
\begin{align}\nonumber
 |\nabla u|^n - \Big(\frac{n-1}{n}\Big)^n\frac{|v|^n}{|x|^n}X_1
  & \geq
   \frac{1}{\kappa_n}
    |x|^{2-n}|v|^{n-2}|\nabla v|^2X_1^{-1}
     \\ \nonumber & \qquad - \Big(\frac{n-1}{n}\Big)^{n-1}
      |x|^{-n}\nabla\big(|v|^n\big)\cdot x.
\end{align}
This means
\begin{align}\nonumber
 I_n[u;\Om]
  & \geq \frac{1}{\kappa_n}
   \int_\Om|x|^{2-n}|v|^{n-2}|\nabla v|^2X_1^{-1}~\dd x
    \\ \nonumber & \qquad + \Big(\frac{n-1}{n}\Big)^{n-1}
     \int_\Om|v|^n\Div\big\{|x|^{-n}x\big\}~\dd x.
 \end{align}
Since $\Div\big\{|x|^{-n}x\big\}=0$ in $\Om\setminus\{0\}$, we deduce \eqref{link2}.
Note that the proof of \eqref{link} in \cite{PsSp} follows the same argument but uses the vectorial inequality
\begin{align}\label{vec old}
 |b-a|^n-|a|^n
  \geq
   \frac{1}{\lambda_n}|b|^n-n|a|^{n-2}a\cdot b,
\end{align}
instead (see also \cite{Li}). \qed
\subsection{An identity for the improved $L^2$-Hardy difference}
One more ingredient we will use is the following equality, originally due to Filippas and Tertikas; see \cite[equality (6.7)]{FT}. We include the proof for the convenience of the reader.

\begin{proposition}\label{Proposition FT 6.7} Let $n\in\mathbb{N}$ and let $U\subset\Rn$ be a bounded domain containing the origin. For any $f\in\test(U\setminus\{0\})$ we have
\begin{align}\nonumber &
 \int_U|\nabla f|^2~\dd x
  -\Big(\frac{n-2}{2}\Big)^2\int_U\frac{|f|^2}{|x|^2}~\dd x
   -\frac14\int_U\frac{|f|^2}{|x|^2}X_1^2~\dd x
 \\ \label{ineq ground state improved} & =
 \int_U|x|^{2-n}|\nabla g|^2X_1^{-1}~\dd x,
\end{align}
where $g$ is defined through $f=|x|^{1-n/2}X_1^{-1/2}g$.
\end{proposition}

\noindent\textbf{Proof.} We compute first
\begin{align}\nonumber
 \nabla f
  =
   -\frac{n-2}{2}|x|^{-n/2}X_1^{-1/2}g\frac{x}{|x|}
    -
     \frac12|x|^{-n/2}X_1^{1/2}g\frac{x}{|x|}
      +
       |x|^{1-n/2}X_1^{-1/2}\nabla g,
\end{align}
in $U\setminus\{0\}$. Expanding now the square, we obtain
\begin{align}\nonumber
 |\nabla f|^2
  & =
   \Big(\frac{n-2}{2}\Big)^2\frac{|f|^2}{|x|^2}
    +
     \frac14\frac{|f|^2}{|x|^2}X_1^2
      +
       |x|^{2-n} X_1^{-1}|\nabla g|^2
        \\ \nonumber & \hspace{1em} +
         \frac{n-2}{2}\frac{g^2}{|x|^n}
          - \frac{n-2}{2}|x|^{1-n}X_1^{-1}\nabla g^2\cdot\frac{x}{|x|}
           - \frac12|x|^{1-n}\nabla g^2\cdot\frac{x}{|x|},
\end{align}
in $U\setminus\{0\}$. But this readily says that the difference between the left hand side and right hand side of \eqref{ineq ground state improved} is given by
\begin{align}\nonumber &
\frac{n-2}{2}\int_U\frac{g^2}{|x|^n}~\dd x
 -\frac12\int_U\big(n-2+X_1\big)|x|^{1-n}X_1^{-1}\nabla g^2\cdot\frac{x}{|x|}~\dd x
\\ \nonumber & \hspace{1em} =
 \frac{n-2}{2}\int_U\frac{g^2}{|x|^n}~\dd x
 +\frac12\int_U\Div\Big\{\big(n-2+X_1\big)|x|^{1-n}X_1^{-1}\frac{x}{|x|}\Big\} g^2~\dd x
\\ \nonumber & \hspace{2em} = 0,
\end{align}
because $\Div\big\{(n-2+X_1)|x|^{1-n}X_1^{-1}\frac{x}{|x|}\big\}=-(n-2)|x|^{-n}$ in $U\setminus\{0\}$. \qed

\section{Proof of the main result}

In this section we give the detailed proof of Theorem \ref{main theorem}. Our proof is reminiscent of the arguments from \cite{GkPs}, that were used to prove the optimal $L^p$-Hardy-Sobolev inequality for $2<p<n$.

\subsection{The case of radial functions} \label{sec:radial}
We start with an elementary lemma which substitutes an end point case of \cite[Lemma 6.1]{Ps} (see also \cite[Lemma 5]{CR}).

\begin{lemma}\label{lemma onedim} For any
$g\in AC\big([0,1]\big)$ with $g(1)=0$ there holds
\begin{align}\label{onedimlemma stepone INEQ}
 \sup_{r\in[0,1]}\Big\{|g(r)|X_2^{1/2}(r)\Big\}
  & \leq
   \Bigg(\int_0^1 t|g'(t)|^2X_1^{-1}(t)~\dd t\Bigg)^{1/2}.
\end{align}
 \end{lemma}

\noindent\textbf{Proof.} Let $0\leq r<1$. Since $f(1)=0$ we have
\begin{align}\nonumber
 |g(r)|
  & \leq
   \int_r^1|g'(t)|~\dd t
    \\ \nonumber & =
     \int_r^1\Big\{t^{-1/2}X_1^{1/2}(t)\Big\}\Big\{t^{1/2}|g'(t)|X_1^{-1/2}(t)\Big\}~\dd t
      \\ \nonumber & \leq
       \Bigg(\int_r^1t^{-1}X_1(t)~\dd t\Bigg)^{1/2}
        \Bigg(\int_r^1t|g'(t)|^2X_1^{-1}(t)~\dd t\Bigg)^{1/2}.
\end{align}
Since $\big(\log X_1(t)\big)' = t^{-1}X_1(t)$, the left integral is easily computed. We deduce
\begin{align}\label{onedimlemma stepone INEQ assist}
 |g(r)|
  \leq
   \Big(-\log X_1(r)\Big)^{1/2}
    \Bigg(\int_0^1t|g'(t)|^2X_1^{-1}(t)~\dd t\Bigg)^{1/2}.
\end{align}
The definition of $X_2$ implies
\[
-\log X_1(r)
  =
   \frac{1-X_2(r)}{X_2(r)}
    \leq
     \frac{1}{X_2(r)},
\]
which when inserted in \eqref{onedimlemma stepone INEQ assist} readily gives \eqref{onedimlemma stepone INEQ}. \qed

\bigskip

The proof of theorem \ref{main theorem} for radial functions is based on the following key proposition.

\begin{proposition}\label{proposition key} Let $q\geq1$. For any radial function $f\in\test(\ub)$ we have
\begin{align}\nonumber &
 \Bigg(\aver_{\hspace{-.3em}\ub}\Big(|f|^{2/n}X_1^{-1-1/n}X_2^{1/n}\Big)^q~\dd x\Bigg)^{1/q}
  \\ \label{radial part} & \leq
   \frac{e^{n/q}}{n}\Gamma^{1/q}\Big(1+q\frac{n-1}{n}\Big)\Bigg(\aver_{\hspace{-.3em}\ub}|x|^{2-n}|\nabla f|^2X_1^{-1}~\dd x\Bigg)^{1/n},
\end{align}
where $\Gamma$ is the gamma function.
\end{proposition}

\noindent\textbf{Proof.} Write $\tilde{f}$ for the radial profile of $f$. From Lemma \ref{lemma onedim} with  $g=\tilde{f}$ we obtain
\begin{align}\nonumber
 |\tilde{f}(r)|^{2/n}X_2^{1/n}(r)
  \leq
   \Bigg(\int_0^1t|\tilde{f}'(t)|^2X_1^{-1}(t)~\dd t\Bigg)^{1/n}.
\end{align}
This is the same as
\begin{align}\nonumber &
 |f(x)|^{2/n}X_2^{1/n}(|x|)
  \leq
   \Bigg(\frac{1}{n\omn}\int_{\ub}|x|^{2-n}|\nabla f|^2X_1^{-1}~\dd x\Bigg)^{1/n}.
\end{align}
Multiplying both sides with $X^{-1+1/n}(|x|)$ and taking $L^q(\ub)$-norms we arrive at

\begin{align}\nonumber &
 \Bigg(\aver_{\hspace{-.3em}\ub}\Big(|f|^{2/n}X_1^{-1+1/n}X_2^{1/n}\Big)^q~\dd x\Bigg)^{1/q}
  \\ \label{progamma} & \leq n^{-1/n}
   \Bigg(\aver_{\hspace{-.3em}\ub}X_1^{-q(1-1/n)}~\dd y\Bigg)^{1/q}\Bigg(\aver_{\hspace{-.3em}\ub}|x|^{2-n}|\nabla f|^2X_1^{-1}~\dd x\Bigg)^{1/n}.
\end{align}
Next we use elementary calculus to estimate the first integral on the right of \eqref{progamma}. Clearly,
\begin{align}\nonumber
 \aver_{\hspace{-.3em}\ub}X_1^{-q(1-1/n)}(|y|)~\dd y
  & =
   n\int_0^1t^{n-1}X_1^{-q(1-1/n)}(t)~\dd t
 \\ \nonumber & =
  \frac{e^n}{n^{q(1-1/n)}}\int_n^\infty e^{-r}r^{q(1-1/n)}~\dd r,
\end{align}
where we have performed the change of variables $t\mapsto e^{1-t/n}$ to reach the last expression. From this we readily get
\begin{align}\label{Gamma}
 \aver_{\hspace{-.3em}\ub}X_1^{-q(1-1/n)}(|y|)~\dd y
  \leq
   \frac{e^n}{n^{q(1-1/n)}}\Gamma\Big(1+q\frac{n-1}{n}\Big).
\end{align}
The proof of \eqref{radial part} follows by inserting \eqref{Gamma} in \eqref{progamma}. \qed

\begin{remark} We note at this point that applying Stirling's formula and elementary estimates, we deduce from \eqref{radial part} that:
\begin{align}\nonumber &
 \Bigg(\aver_{\hspace{-.3em}\ub}\Big(|f|^{2/n}X_1^{-1-1/n}X_2^{1/n}\Big)^q~\dd x\Bigg)^{1/q}
  \\ \label{radial part sobolev growth} & \leq
   c(n)~q^{1-1/n}\Bigg(\aver_{\hspace{-.3em}\ub}|x|^{2-n}|\nabla f|^2X_1^{-1}~\dd x\Bigg)^{1/n},
\end{align}
for any radial function $f\in\test(\ub)$. This is a form of \eqref{radial part} in the spirit of \eqref{sharp Trudinger growth}. In fact, it is \eqref{radial part sobolev growth} we are utilizing in the next section to treat the spherical mean of a general function $f\in\test(\ub)$ (not necessarily radial).
\end{remark}

\bigskip

\noindent\textbf{Proof of theorem \ref{main theorem} for radial functions.} Let $u\in\mathcal{I}_n(\ub)\setminus\{0\}$ be radially symmetric and define $v$ through $u=vX_1^{-1+1/n}$. The function
\[
 w:=|v|^{n/2},
\]
is then radially symmetric and also $w\in\test(\ub)$. Therefore, by taking $f=w$ in proposition \ref{proposition key} we obtain
\begin{align}\nonumber &
 \Bigg(\aver_{\hspace{-.3em}\ub}\Big(|v|X_1^{-1-1/n}X_2^{1/n}\Big)^q~\dd x\Bigg)^{1/q}
  \\ \nonumber & \leq \Big(\frac n2\Big)^{2/n}
   \frac{e^{n/q}}{n}\Gamma^{1/q}\Big(1+q\frac{n-1}{n}\Big)
    \Bigg(\aver_{\hspace{-.3em}\ub}|x|^{2-n}|v|^{n-2}|\nabla v|^2X_1^{-1}~\dd x\Bigg)^{1/n}.
\end{align}
Taking into account \eqref{link2} with the hypothesis $I_n[u;\ub]\leq1$, this becomes
\begin{align}\label{L^q presum}
 \Bigg(\aver_{\hspace{-.3em}\ub}\Big(|u|X_2^{1/n}\Big)^q~\dd x\Bigg)^{1/q}
  \leq
   e^{n/q}\bigg(\frac{\kappa_n}{4\omn n^{n-2}}\bigg)^{1/n}
    \Gamma^{1/q}\Big(1+q\frac{n-1}{n}\Big).
\end{align}
The rest of the proof is standard: Taking $q=ln/(n-1)$, $l\in\N$, in \eqref{L^q presum} we get
\begin{align}\nonumber
 \aver_{\hspace{-.3em}\ub}\Big[\Big(|u(x)|X_2^{1/n}(|x|)\Big)^{n/(n-1)}\Big]^l~\dd x
 \leq
   e^n\bigg(\frac{\kappa_n}{4\omn n^{n-2}}\bigg)^{l/(n-1)}\Gamma(1+l).
\end{align}
Note that this holds true also for $l=0$. Now we multiply both sides by $a^l/l!$, note that $\Gamma(1+l)=l!$ and add for all integers $l\in[0,m]$, $m\in\N$, to arrive at
\begin{align}\nonumber
 \aver_{\hspace{-.3em}\ub}\sum_{l=0}^m
  \frac{1}{l!}\Big[a\Big(|u(x)|X_2^{1/n}(|x|)\Big)^{n/(n-1)}\Big]^l\dd x
   \leq
    e^n\sum_{l=0}^m\Bigg[a
     \bigg(\frac{\kappa_n}{4\omn n^{n-2}}\bigg)^{1/(n-1)}
      \Bigg]^l.
\end{align}
The series on the right converges if and only if
\begin{align}\label{estimalpha}
 a<\bigg(\frac{4\omn n^{n-2}}{\kappa_n}\bigg)^{1/(n-1)}.
\end{align}
Hence for any such $a$ the proof is completed by letting $m\rightarrow\infty$ and using the monotone convergence theorem. \qed

\smallskip

\begin{remark} Since $\kappa_2=1$, for $n=2$, estimate \eqref{estimalpha} reads $a<4\pi$; that is $a$ has to be strictly less than Moser's sharp constant. We don't know if it still holds true for $a=4\pi$.
\end{remark}
\subsection{Proof of Theorem \ref{main theorem}} \label{sec:non-radial}
 Hence it is enough to establish \eqref{sharp Leray-Trudinger growth} for $\Om=B_{R_\Om}(0)$. Furthermore, being scaling invariant, it is enough to consider only the case $R_\Om=1$. Finally, it is enough to assume that $u\in\test(\ub\setminus\{0\})$ (since $W_0^{1,n}(\ub\setminus\{0\})=W_0^{1,n}(\ub)$; see for instance \cite[Theorem 2.43]{HKM}). Pick $u\in\mathcal{I}_n(\ub\setminus\{0\})\setminus\{0\}$ and consider the function $v$ defined through the transformation $u=vX_1^{-1+1/n}$. Following \cite{VzZ}, we use spherical coordinates $x=(r,\theta)$ ($r=|x|$ and $\theta=x/|x|$) to decompose $v$ into spherical harmonics. For this purpose, let $\{h_l\}_{l\in\N\cup\{0\}}$ be the orthonormal basis of $L^2(\us)$ that is comprised of eigenfunctions of the Laplace-Beltrami operator $-\Delta_{\us}$ (the angular part of the Laplacian when expressed in spherical coordinates). This has corresponding eigenvalues $\lambda_l=l(l+n-2)$, $l\in\N\cup\{0\}$ (see \cite[Appendix]{Schn}). Thus for all $l,m\in\N\cup\{0\}$ we have
\[
 -\Delta_{\us}h_l=\lambda_lh_l~~\mbox{on }\us
  \quad\mbox{ and }\quad
   \aver_{\hspace{-.3em}\us}h_l(\theta)h_m(\theta)~\dS(\theta)=\delta_{lm}.
\]
With these definitions we have the decomposition of $v$ in its spherical harmonics
\[
 v(x)=\sum_{l=0}^\infty v_l(r)h_l(\theta).
\]
In particular $h_0(\theta)=1$ and the first term in the above decomposition is given by the spherical mean of $v$ on $\prt B_r(0)$, that is
\[
 v_0(r)=\aver_{\hspace{-.3em}\prt B_r(0)}v(x)~\dS(x)=\aver_{\hspace{-.3em}\us}v(r\theta)~\dS(\theta).
\]
Now let $q>n$. Minkowski's inequality implies
\begin{align}\nonumber &
 \Bigg(\aver_{\hspace{-.3em}\ub}\Big(|u|X_2^{1/n}\Big)^q~\dd x\Bigg)^{1/q}
  \\ \label{Mink} & \leq
   \Bigg(\aver_{\hspace{-.3em}\ub}\Big(|u_0|X_2^{1/n}\Big)^q~\dd x\Bigg)^{1/q}
   +
   \Bigg(\aver_{\hspace{-.3em}\ub}\Big(|u-u_0|X_2^{1/n}\Big)^q~\dd x\Bigg)^{1/q},
\end{align}
where $u_0:=v_0X_1^{-1+1/n}$. Note that since $X_1$ is a radial function, $u_0$ is just the spherical mean of $u$ on $\prt B_r(0)$; a radial function.

\medskip

Define $w:=|v|^{n/2}$ and consider $w_0$; that is, the spherical mean of $w$ on $\prt B_r(0)$. H\"{o}lder's inequality implies $|v_0|\leq w_0^{2/n}$. Indeed,
\begin{align}\nonumber
 |v_0(r)|^{n/2}
  & = \bigg|\aver_{\hspace{-.3em}\prt B_r(0)}v(y)~\dS(y)\bigg|^{n/2}
   \\ \nonumber & \leq
    \aver_{\hspace{-.3em}\partial B_r(0)}|v(y)|^{n/2}~\dS(y)
     \\ \nonumber & =
      \aver_{\hspace{-.3em}\prt B_r(0)}w(y)~\dS(y)
       = w_0(r).
\end{align}
Using this, we see for the first term in \eqref{Mink} that
\begin{align}\nonumber
 \Bigg(\aver_{\hspace{-.3em}\ub}
 \Big(|u_0|X_2^{1/n}\Big)^q~\dd x\Bigg)^{1/q}
  = &
 \Bigg(\aver_{\hspace{-.3em}\ub}
 \Big(|v_0|X_1^{-1+1/n}X_2^{1/n}\Big)^q~\dd x\Bigg)^{1/q}
  \\ \nonumber & \hspace{-2em} \leq
 \Bigg(\aver_{\hspace{-.3em}\ub}
 \Big(|w_0|^{2/n}X_1^{-1+1/n}X_2^{1/n}\Big)^q~\dd x\Bigg)^{1/q}
  \\ \label{Mink LEFT} & \hspace{-3em} \leq
   c_1(n)~q^{1-1/n}
   \Bigg(\int_\ub|x|^{2-n}|\nabla w_0|^2X_1^{-1}~\dd x\Bigg)^{1/n},
\end{align}
because of \eqref{radial part sobolev growth} with $f=w_0$.

For the second term of \eqref{Mink} we use first the fact that $X_2\leq1$ and then \eqref{sharp Trudinger growth} to the function $u-u_0\in W_0^{1,n}(\ub)$ to get
\begin{align}\nonumber&
 \Bigg(\aver_{\hspace{-.3em}\ub}\Big(|u-u_0|X_2^{1/n}\Big)^q~\dd x\Bigg)^{1/q}
  \\ \nonumber & \leq
   c_2(n)~q^{1-1/n}
   \Bigg(\int_{\ub}\big|\nabla(u-u_0)\big|^n~\dd x\Bigg)^{1/n}
    \\ \label{Mink RIGHT} & =
      c_2(n)~q^{1-1/n}
     \Bigg(\int_{\ub}\big|\nabla\big[(v-v_0)X_1^{-1+1/n}\big]\big|^n~\dd x\Bigg)^{1/n}.
\end{align}
\textbf{Claim:} The assumption $I_n[u]\leq1$ is enough for each of the two integrals on the right of \eqref{Mink LEFT} and \eqref{Mink RIGHT} to be bounded by a constant not depending on $u$. In particular we will show that
\begin{align}\label{Step1}
 \int_\ub|x|^{2-n}|\nabla w_0|^2X_1^{-1}~\dd x
  \leq c_3(n)
   \int_\ub|x|^{2-n}|v|^{n-2}|\nabla v|^2X_1^{-1}~\dd x,
\end{align}
and that
\begin{align}\label{Step2}
 \int_{\ub}\big|\nabla\big[(v-v_0)X_1^{-1+1/n}\big]\big|^n~\dd x
  \leq c_4(n)
   \int_\ub|\nabla v|^nX_1^{-n+1}~\dd x.
\end{align}
Because of \eqref{link} and \eqref{link2}, estimates \eqref{Step1} and \eqref{Step2} will readily imply our claim.

\medskip

\noindent\textbf{Proof of \eqref{Step1}:}
Set $\zeta:=|x|^{1-n/2}X_1^{-1/2}w$, then by proposition \ref{Proposition FT 6.7} with $U=\ub$ and $f=\zeta_0$, we have the following equality
\begin{align}\nonumber &
 \int_\ub|x|^{2-n}|\nabla w_0|^2X_1^{-1}~\dd x
  \\ \label{w0 to improved0} & =
  \int_\ub|\nabla \zeta_0|^2~\dd x
   -\Big(\frac{n-2}{2}\Big)^2\int_\ub\frac{|\zeta_0|^2}{|x|^2}~\dd x
    -\frac14\int_\ub\frac{|\zeta_0|^2}{|x|^2}X_1^2~\dd x.
\end{align}
From \cite[eq. (7.6)]{FT}) we know that
\begin{align}\nonumber &
 \int_\ub|\nabla \zeta_0|^2~\dd x
  -\Big(\frac{n-2}{2}\Big)^2\int_\ub\frac{|\zeta_0|^2}{|x|^2}~\dd x
   -\frac14\int_\ub\frac{|\zeta_0|^2}{|x|^2}X_1^2~\dd x
 \\ \nonumber & \leq \int_\ub|\nabla \zeta|^2\dd x
   -\Big(\frac{n-2}{2}\Big)^2\int_\ub\frac{|\zeta|^2}{|x|^2}~\dd x
    -\frac14\int_\ub\frac{|\zeta|^2}{|x|^2}X_1^2~\dd x
 \\ \label{improved0 to improved and w} & =
  \int_\ub|x|^{2-n}|\nabla w|^2X_1^{-1}~\dd x,
\end{align}
the last equality again because of proposition \ref{Proposition FT 6.7} with $U=\ub$, but $f=\zeta$ this time. We have just showed that
\begin{align}\nonumber
 \int_\ub|x|^{2-n}|\nabla w_0|^2X_1^{-1}~\dd x
  \leq
   \int_\ub|x|^{2-n}|\nabla w|^2X_1^{-1}~\dd x,
\end{align}
and the proof of \eqref{Step1} follows at once since $w:=|v|^{n/2}$.

\medskip

\noindent\textbf{Proof of \eqref{Step2}:} We start from the right hand side of \eqref{Step2}, by noticing that
\begin{align}\nonumber
 & \int_\ub|\nabla v|^nX_1^{-n+1}~\dd x
  \\ \nonumber & \quad = \int_0^1X_1^{-n+1}(r)~r^{n-1}
  \int_{\us}\bigg((\prt_rv)^2+\frac{1}{r^2}|\nabla_\theta v|^2\bigg)^{n/2}
  ~\dS(\theta)~\dd r
   \\ \nonumber & \quad \geq
    \int_0^1X_1^{-n+1}(r)~r^{n-1}\int_{\us}|\prt_rv|^n~\dS(\theta)~\dd r
     \\ \label{A1A2} & \qquad +\underbrace{\int_0^1X_1^{-n+1}(r)~r^{-1}
     \int_{\us}|\nabla_\theta v|^n
     ~\dS(\theta)~\dd r}_{=:J},
\end{align}
by the fact that $(\kappa+\lambda)^q\geq \kappa^q+\lambda^q$, for all $\kappa,\lambda\geq0$ and any $q\geq1$, and also because $X_1^{-n+1}(r)\geq1$ for all $r$. To estimate the first term on the right of \eqref{A1A2} we use \eqref{vec old} to get
\begin{align}\nonumber
 \int_{\us}&|\prt_rv|^n~\dS(\theta)
  \\ \nonumber & \geq \int_{\us}|\prt_rv_0|^n~\dS(\theta) + \frac{1}{\lambda_n}\int_{\us}|\prt_r(v-v_0)|^n~\dS(\theta)
   \\\label{A11} & \quad + n\int_{\us}|\prt_rv_0|^{n-2}(\prt_r v_0)\prt_r(v-v_0)~\dS(\theta).
\end{align}
But since $\{v_l\}_{l\in\N\cup\{0\}}$ are radial
\begin{align*}&
\int_{\us}|\prt_rv_0|^{n-2}(\prt_r v_0)\prt_r(v-v_0)~\dS(\theta)
 \\ & \qquad = |v'_0(r)|^{n-2}v'_0(r)\int_{\us}\partial_r(v-v_0)~\dS(\theta)
  \\ & \qquad = |v'_0(r)|^{n-2}v'_0(r)\sum_{l=1}^\infty v'_l(r)\int_{\us}h_l(\theta)~\dS(\theta)
    = 0,
\end{align*}
and so
\[
 \int_{\us}|\prt_rv|^n~\dS(\theta)
  \geq
   \frac{1}{\lambda_n}\int_{\us}|\prt_r(v-v_0)|^n~\dS(\theta),
\]
where we have cancel also the first term on the right hand side of \eqref{A11}. Plugging this to \eqref{A1A2} we deduce
\begin{align}\nonumber &
 \int_\ub|\nabla v|^nX_1^{-n+1}~\dd x
  \\ \nonumber & \geq
   \frac{1}{\lambda_n}\int_0^1X_1^{-n+1}(r)~r^{n-1}
   \int_\us\bigg(|\prt_r(v-v_0)|^n+\frac{1}{r^n}|\nabla_\theta v|^n\bigg)
   ~\dS(\theta)~\dd r
 \\ \nonumber & \qquad + \Big(1-\frac{1}{\lambda_n}\Big)J
 \\ \nonumber & \geq \frac{2^{1-n/2}}{\lambda_n} \int_0^1X_1^{-n+1}(r)~r^{n-1}\int_{\us}\bigg(|\prt_r(v-v_0)|^2+\frac{1}{r^2}|\nabla_\theta v|^2\bigg)^{n/2}
  ~\dS(\theta)~\dd r
 \\ \nonumber & \qquad +\Big(1-\frac{1}{\lambda_n}\Big)J
  \\ \label{A1A2A3} & =
   \frac{2^{1-n/2}}{\lambda_n}\int_\ub|\nabla(v-v_0)|^nX_1^{-n+1}~\dd x
   +\Big(1-\frac{1}{\lambda_n}\Big)J,
\end{align}
by the fact that $\kappa^q+\lambda^q\geq2^{1-q}(\kappa+\lambda)^q$, for all $\kappa,\lambda\geq0$ and any $q\geq1$. To estimate $J$ observe first that

\[
 \int_{\us}(v-v_0)~\dS(\theta)=\sum_{l=1}^\infty v_l(r)\int_{\us}h_l(\theta)~\dS(\theta)=0,
\]
so that utilizing once more the fact that $v_0$ is radial, we may use the Poincar\'e inequality on $\us$ (see for example \cite{H})
\begin{align}\nonumber
 \int_{\us}|\nabla_\theta v|^n~\dS(\theta)
  & = \int_{\us}|\nabla_\theta(v-v_0)|^n~\dS(\theta)
  \\ \nonumber & \geq
   C_{\mathbf{P}}(n)\int_{\us}|v-v_0|^n~\dS(\theta).
\end{align}
Inserting this in the definition of $J$, we get from \eqref{A1A2A3} the existence of a positive constant $C(n)$ such that
\begin{align}\nonumber &
 \int_\ub|\nabla v|^nX_1^{-n+1}~\dd x
  \\ \nonumber & \geq C(n)\Bigg(\int_\ub|\nabla (v-v_0)|^nX_1^{-n+1}~\dd x
   +
    \int_\ub|x|^{-n}|v-v_0|^nX_1^{-n+1}~\dd x\Bigg)
\end{align}
Since $X_1\leq1$ we replace $X_1^{-n+1}$ by $X_1$ in the second integral to deduce
\begin{align}\nonumber &
 \int_\ub|\nabla v|^nX_1^{-n+1}~\dd x
  \\ \nonumber & \geq C(n)\Bigg(\int_\ub|\nabla (v-v_0)|^nX_1^{-n+1}~\dd x
   +
    \int_\ub|x|^{-n}|v-v_0|^nX_1~\dd x\Bigg)
     \\ \nonumber & \geq C'(n)\int_\ub\big|\nabla\big[(v-v_0)X_1^{-1+1/n}\big]\big|^n~\dd x
\end{align}
as required. \qed
\begin{acknowledgements}
G. di Blasio and G. Pisante are members of the Gruppo Nazionale per l'Analisi Matematica, la Probabilit\`{a} e le loro Applicazioni (GNAMPA) of the Istituto Nazionale di Alta Matematica (INdAM) whose support is gratefully acknowledged. The research has also been supported by project Vain-Hopes within the program VALERE: VAnviteLli pEr la RicErca. 
\end{acknowledgements}

\end{document}